\documentclass{ifacconf}
\usepackage{color}
\usepackage{xurl}

\usepackage{amsmath}

\usepackage{dsfont}
\usepackage[]{mdframed}
\usepackage{amssymb}
\usepackage{algorithm} 
\usepackage{algpseudocode} 

\usepackage{graphicx}      
\usepackage{natbib}        

\allowdisplaybreaks

\begin{document}
\begin{frontmatter}


\title{REAP-T: A MATLAB Toolbox for Implementing Robust-to-Early Termination Model Predictive Control}

\author{Mohsen Amiri and Mehdi Hosseinzadeh}

\address{School of Mechanical and Materials Engineering, Washington State University, Pullman, WA 99164, USA \\
E-mail: mohsen.amiri@wsu.edu;~mehdi.hosseinzadeh@wsu.edu.}

\begin{abstract}  
This paper presents a MATLAB toolbox for implementing robust-to-early termination model predictive control, abbreviated as REAP, which is designed to ensure a sub-optimal yet feasible solution when MPC computations are prematurely terminated due to limited computational resources. Named REAP-T, this toolbox is a comprehensive, user-friendly, and modular platform that enables users to explore, analyze, and customize various components of REAP for their specific applications. Notable attributes of REAP-T are: (i) utilization of built-in MATLAB functions for defining the MPC problem; (ii) an interactive and intuitive graphical user interface for parameter tuning and visualization; (iii) real-time simulation capabilities, allowing users to observe and understand the real-time behavior of their systems; and (iv) inclusion of real-world examples designed to guide users through its effective use. 
\end{abstract}


\end{frontmatter}

\section{Introduction}\vspace{-0.5cm}
Developing practical control strategies for real-world systems with constraints has become a significant challenge in recent years. Model Predictive Control (MPC) \citep{rawlings2017model, camacho2007model} is a widely adopted approach that addresses this challenge by optimizing performance objectives over a receding horizon at each sampling instant, while ensuring the satisfaction of all safety and operational constraints. However, despite its effectiveness, MPC's reliance on online optimization introduces significant implementation challenges, especially for systems with fast dynamics or limited processing power.

Over the past few decades, various approaches have been proposed to address the computational challenges associated with MPC. A notable recent advancement is the Robust-to-Early terminAtion oPtimization (REAP) theory introduced in \citep{Hosseinzadeh2023RobustTermination}. REAP embeds the MPC solution within a virtual continuous-time dynamical system based on primal-dual gradient flow \citep{Feijer2010,HosseinzadehTAC2020,Hosseinzadeh2022_ROTEC,LucaLetter,HosseinzadehLetter,HosseinzadehTAC2024}, ensuring that even when execution is prematurely terminated, the resulting solution remains suboptimal yet feasible. This adaptability allows REAP to effectively handle limited and varying computational resources, while maintaining feasibility and achieving control objectives. The REAP theory was further developed in \citep{amiri2025practical} to facilitate its implementation across various computing platforms. The study provided implementation details for two illustrative examples, along with the corresponding MATLAB files.




We believe that REAP has reached a level of maturity that positions it as a standard tool for implementing MPC in systems with limited computational resources. To support its adoption, this paper introduces a comprehensive MATLAB toolbox, called REAP-T, designed to facilitate the application of REAP to constrained control problems. The key features of REAP-T are as follows: (1) It computes control inputs in real time, adapting to the available computation time on the hardware at each sampling instant; (2) it handles linear systems with linear constraints on states and inputs, allowing users to apply it to their specific systems; (3) it accommodates a variety of control objectives, including output tracking and convergence to equilibrium points; and (4) it allows users to select the strategy for developing the terminal constraint set and to investigate how their selection impacts closed-loop performance.


The remainder of this paper is organized as follows: Section \ref{sec:REAP theory} provides a summary of the REAP theory. Section \ref{sec:REAP-Toolbox} introduces REAP-T and its functionalities. Section \ref{sec:Tutorial} demonstrates the toolbox with two real-world examples. Finally, Section \ref{sec:Conclusion} concludes the paper.

\section{REAP Theory}\label{sec:REAP theory}

Consider the following continuous-time Linear Time-Invariant (LTI) dynamic system:
\begin{align}\label{eq:SystemContinous}
\dot{x}(t) = A_c x(t) + B_c u(t),~~~~~~y(t) = C_c x(t) + D_c u(t),
\end{align}
where $x(t)=[x_1(t)~\cdots~x_n(t)]^\top \in \mathbb{R}^n$ represents the state vector, $u(t)=[u_z(t)~\cdots~u_p(t)]^\top \in \mathbb{R}^p$ is the control input, $y(t)=[y_1(t)~\cdots~y_m(t)]^\top\in \mathbb{R}^m$ is the output vector, $t$ is the continuous time variable, and $A_c \in \mathbb{R}^{n \times n}$, $B_c \in \mathbb{R}^{n \times p}$, $C_c \in \mathbb{R}^{m \times n}$, and $D_c \in \mathbb{R}^{m \times p}$ are system matrices.

Although the control input $u(t)$ should be applied continuously to the system, its computation is usually performed in discrete time. For a zero-order hold implementation with a sampling period of $\Delta T\in\mathbb{R}_{>0}$ seconds, the discrete-time model can be expressed as:
\begin{align}\label{eq:System1}
x(k+1) &= A x(k) + B u(k),~~y(k)= C x(k) + D u(k),
\end{align}
where $k$ denotes the sampling instant (i.e., $[k,k+1)$ is equal to $\Delta T$ seconds), $A=e^{A_c\Delta T}$, $B=\int_{0}^{\Delta T}e^{A_ct}B_cdt$, $C=C_c$, and $D=D_c$. We assume that ($A, B$) is controllable.  


The states and inputs of system \eqref{eq:System1} are subject to the following constraints at all times:
\begin{align}\label{eq:constraints}
x(k) \in \mathcal{X} \subset \mathbb{R}^n,~~~~~u(k)\in \mathcal{U} \subset \mathbb{R}^p,~~~~~\forall k\in\mathbb{Z}_{\geq0},
\end{align}
where $\mathcal{X}$ and $\mathcal{U}$ are convex polytopes defined as:
\begin{subequations}\label{eq:allConstraints}
\begin{align}
\mathcal{X} &= \left\{ x \in \mathbb{R}^n : \underline{x}_i\leq x_i \leq \bar{x}_i,~i = 1, \ldots, n \right\}, \label{eq:Xconstraints} \\
\mathcal{U} &= \left\{ u \in \mathbb{R}^p : \underline{u}_i\leq u_i \leq \bar{u}_i,~i = 1, \ldots, p \right\}, \label{eq:Uconstraints}
\end{align}
\end{subequations}
where $\underline{x}_i\in \mathbb{R}$ and $\bar{x}_i \in \mathbb{R}$ represent the lower and upper limits for the $i$-th state, respectively, and $\underline{u}_i\in \mathbb{R}$ and $\bar{u}_i \in \mathbb{R}$ denote the lower and upper limits for the $i$-th control input. 


Let $ r \in \mathbb{R}^m $ denote the desired reference. There exists at least one steady-state configuration $\big(\bar{\mathbf{x}}_r,\bar{\mathbf{u}}_r\big)$ such that:
\begin{align}\label{eq:SteadyState}
 \bar{\mathbf{x}}_r = A \bar{\mathbf{x}}_r + B \bar{\mathbf{u}}_r,~~~~~~r = C \bar{\mathbf{x}}_r + D \bar{\mathbf{u}}_r,
\end{align} 
where $ \bar{\mathbf{x}}_r \in \operatorname{Int}(\mathcal{X}) $ and $ \bar{\mathbf{u}}_r \in \operatorname{Int}(\mathcal{U}) $. This reference signal is termed a steady-state admissible reference, and the set of all such references is denoted by $\mathcal{R}\subset\mathbb{R}^m$.

\subsection{MPC Formulation}
Given the desired reference $r\in\mathcal{R}$ and prediction horizon length $N\in\mathbb{Z}_{>0}$, MPC computes the optimal control sequence $\mathbf{u}^\ast(k) := [u^\ast(0|k)^\top, \ldots, u^\ast(N-1|k)^\top]^\top \in \mathbb{R}^{Np}$ at any time instant $k$ by solving the following optimization problem:
\begin{subequations}\label{eq:OptimizationProblemMain}
\begin{align}
\arg\min_{\mathbf{u}} & \sum_{s=0}^{N-1} \| \hat{x}(s|k) - \bar{\mathbf{x}}_r \|_{Q_x}^2 + \sum_{s=0}^{N-1} \| u(s|k) - \bar{\mathbf{u}}_r \|_{Q_u}^2 \nonumber \\
& + \| \hat{x}(N|k) - \bar{\mathbf{x}}_r \|_{Q_N}^2, \label{eq:CostFunctionExtend}
\end{align}
subject to:
\begin{align}
& \hat{x}(s+1|k) = A \hat{x}(s|k) + B u(s|k), \quad \hat{x}(0|k) = x(k), \label{eq:ConstraintRevised} \\
& \hat{x}(s|k) \in \mathcal{X}, \; s = 0, \ldots,N-1, \\
& u(s|k) \in \mathcal{U}, \; s = 0, \ldots, N-1,  \\
& (\hat{x}(N|k), r) \in \Omega, \label{eq:ConstraintTerminalExtend}
\end{align}
\end{subequations}
where $Q_x \in \mathbb{R}^{n \times n}$, $Q_u \in \mathbb{R}^{p \times p}$, and $Q_N \in \mathbb{R}^{n \times n}$ are positive semi-definite weighting matrices, and $\Omega$ is the terminal constraint set. The computation of the terminal cost matrix $Q_N$ and the terminal constraint set $\Omega$ will be detailed in Section \ref{sec:TerminalConstraintSet}.

\subsection{Robust-to-Early Termination MPC}

At any time instant $k$, REAP uses the primal-dual gradient flow to inform the following virtual discrete-time dynamical system \citep{amiri2025practical}:
\begin{subequations}\label{eq:DisREAP}
\begin{align}
\hat{\mathbf{u}}(\tau|k) =& \hat{\mathbf{u}}(\tau-1|k)-\sigma(\tau|k)\cdot d\tau\cdot \nabla_{\hat{\mathbf{u}}} \mathcal{B}\big(x(k), r, \nonumber\\
&\hat{\mathbf{u}}(\tau-1|k),\hat{\lambda}(\tau-1|k)\big), \label{eq:DisREAP1}
\\
\hat{\lambda}(\tau|k) =& \hat{\lambda}(\tau-1|k)+\sigma(\tau|k)\cdot d\tau\cdot  \Big(\nabla_{\hat{\lambda}} \mathcal{B}\big(x(k), r, \nonumber\\
&\hat{\mathbf{u}}(\tau-1|k),\hat{\lambda}(\tau-1|k)\big)+ \Phi(\tau-1|k) \Big),\label{eq:DisREAP2}
\end{align}
\end{subequations}
where $\tau$ represents the computation step at each time instant, $d\tau$ is the discretization step, $\sigma(\tau|k)$ is the Karush-Kuhn-Tucker (KKT) parameter at computation step $\tau$, $\mathcal{B}\big(x(k),r,\mathbf{u},\lambda\big)$ is the modified barrier function \citep{Polyak1992,Melman1996,Vassiliadis1998} associated with the optimization problem \eqref{eq:OptimizationProblemMain} and $\Phi(\tau|k)$ is the projection operator. It has been shown in \citep{amiri2025practical} that if $\sigma(\tau|k)$ is determined using an adaptive law (Equation (16) of \citep{amiri2025practical}), then $\big(\hat{\mathbf{u}}(\tau|k),\hat{\lambda}(\tau|k)\big)$ converges to $\big(\mathbf{u}^\ast(k),\lambda^{\ast}(k)\big)$ as $\tau \rightarrow \infty$, while ensuring that $\hat{\mathbf{u}}(\tau|k)$ remains a feasible solution for the MPC problem \eqref{eq:OptimizationProblemMain} and $\hat{\lambda}(\tau|k)\in\mathbb{R}_{\geq0}$ at all times $\tau$, where $\mathbf{u}^\ast(k)$ is as in \eqref{eq:OptimizationProblemMain} and $\lambda^{\ast}(k)$ is the optimal dual variable at time instant $k$. 


\section{REAP-T}\label{sec:REAP-Toolbox}

\begin{figure*}[t]
    \centering
     \includegraphics[width=.75\linewidth]{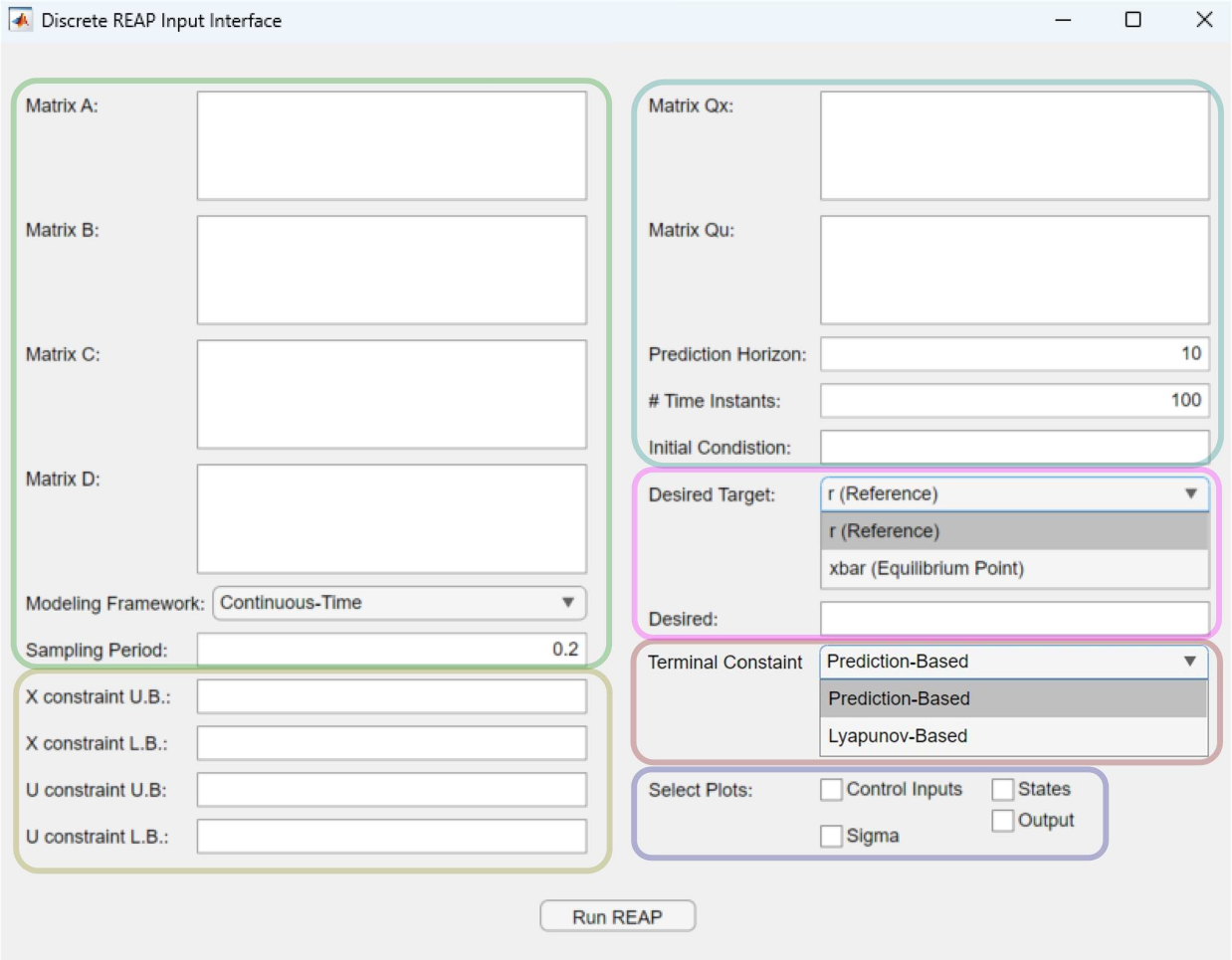}
    \caption{The developed GUI for running REAP to solve the MPC problem.}
    \label{fig:GUI}
    \vspace{0.1cm}
    \hrule
\end{figure*}

This section presents a step-by-step guide to utilizing the functionalities of REAP-T. The modular design of REAP-T encompasses different layers of functions, enabling users to explore modules for in-depth investigation or to leverage the functions for developing their own applications.
 
 
\subsection{Installation}
To install REAP-T, begin by downloading its source files from: \url{https://github.com/mhsnar/REAP-T}. REAP-T incorporates the YALMIP toolbox \citep{Lofberg2004}, which is used exclusively during the design phase to compute the initial feasible solution at time instant $k=0$. After extracting the YALMIP zip file to your desired location, run the \texttt{REAPT\_GUI} function in MATLAB to launch the Graphical User Interface (GUI).

\subsection{System Dynamics\textemdash Green Box in Fig. \ref{fig:GUI}}
The first step in using REAP-T is to define the system dynamics, which includes specifying the system matrices and the sampling period $\Delta T$. sers can choose to define the system in either continuous-time (by specifying matrices $A_c$, $B_c$, $C_c$, and $D_c$) or discrete-time (by specifying matrices $A$, $B$, $C$, and $D$) using a dropdown menu. It is noteworthy that if the system is defined in the continuous domain, REAP-T will utilize the specified sampling period and the zero-order hold method to discretize the system. 

\begin{rem}
If ($A,B$) is not controllable, REAP-T will terminate execution and display the following message:\vspace{0.1cm}
\begin{mdframed}
\texttt{The pair (A,B) is not controllable. REAP-T\\ 
cannot proceed with the specified system.}
\end{mdframed}   
\end{rem}





\subsection{Constraints\textemdash Yellow Box in Fig. \ref{fig:GUI}}

Once the system's dynamics have been defined, constraints on the states and inputs can be defined using the GUI.

\noindent\textbf{Constraints on States:}
State constraints, as defined in \eqref{eq:constraints}, can be specified by setting upper and lower bounds as follows:
\begin{subequations}\label{eq:XConstraintsGUI}
 \begin{align}
\text{X constraint U.B.}=&[\bar{x}_1~\bar{x}_2~\cdots~\bar{x}_n]^\top,\\
\text{X constraint L.B.}=&[\underline{x}_1~\underline{x}_2~\cdots~\underline{x}_n]^\top,
\end{align}   
\end{subequations}
where $\bar{x}_i$ and $\underline{x}_i$ represent the upper and lower bounds for the $i$-th state, respectively. Note that if the $i$-th state is not bounded on either side, the corresponding limit should be set to \texttt{Inf} or \texttt{-Inf} in \eqref{eq:XConstraintsGUI}.


\noindent\textbf{Constraints on Inputs:}
Similarly, input constraints, as defined in \eqref{eq:constraints}, can be specified by setting upper and lower bounds as follows:
\begin{subequations}\label{eq:UConstraintsGUI}
\begin{align}
\text{U constraint U.B.}=&[\bar{u}_1~\bar{u}_2~\cdots~\bar{u}_p]^\top,\\
\text{U constraint L.B.}=&[\underline{u}_1~\underline{u}_2~\cdots~\underline{u}_p]^\top,
\end{align}
\end{subequations}
where $\bar{u}_i$ and $\underline{u}_i$ represent the upper and lower bounds for the $i$-th control input, respectively. If the $i$-th control input is not bounded on either side, the corresponding limit should be set to \texttt{Inf} or \texttt{-Inf} in \eqref{eq:UConstraintsGUI}.

\noindent\textbf{List of Constraints:} Users can verify the defined constraints by using the following command:
\begin{verbatim}
Functions.list(AllConstraints)
\end{verbatim} 

This command will display all constraints on both states and control inputs. For instance:\vspace{0.1cm}
\begin{mdframed}
\texttt{State Constraints:}\\
\texttt{State Constraint 1: x1 $<$= 5}\\
\texttt{State Constraint 2: x2 $<$= Inf}\\
\texttt{State Constraint 3: x1 $>$= -Inf}\\
\texttt{State Constraint 4: x2 $>$= -5}\\
\texttt{Input Constraints:}\\
\texttt{Input Constraint 1: u1 $<$= 10}\\
\texttt{Input Constraint 2: u2 $<$= 10}   
\end{mdframed}




\subsection{Simulation Parameters\textemdash Blue Box in Fig. \ref{fig:GUI}}
For the simulation, the user should specify the simulation duration (denoted by \texttt{\# Time Instants}), the length of the prediction horizon $N$, and weighting matrices $Q_x$ and $Q_u$ as defined in \eqref{eq:CostFunctionExtend} to balance tracking performance and control effort. Additionally, the user needs to specify the initial condition $x(0)\in\mathbb{R}^n$. If the given initial condition does not lie within the region of attraction of the MPC problem, REAP-T displays the following message:
\begin{mdframed}
\texttt{The specified initial condition does not \\ 
belong to the region of attraction.\\
REAP-T cannot proceed.}
\end{mdframed}

Note that to obtain the terminal cost matrix $Q_N$, REAP-T uses the following Riccati equation: $Q_N = A^\top Q_N A - (A^\top Q_N B)(Q_u + B^\top Q_N B)^{-1}(B^\top Q_N A) + Q_x$. As discussed in \citep{nicotra2018embedding,Hosseinzadeh2023RobustTermination}, this selection is optimal for the unconstrained problem.

\subsection{Desired Target\textemdash Pink Box in Fig. \ref{fig:GUI}}
Users can choose the desired target type\textemdash either a desired reference $r$ or desired equilibrium point $\bar{\mathbf{x}}_r$\textemdash  from a dropdown menu. If the desired reference $r$ is specified, the function \texttt{desiredCalculation} computes the steady-state configuration $\big(\bar{\mathbf{x}}_r,\bar{\mathbf{u}}_r\big)$ by solving the equations in \eqref{eq:SteadyState}. Alternatively, if $\bar{\mathbf{x}}_r$ is provided, REAP-T solves \eqref{eq:SteadyState} to determine the corresponding steady input $\bar{\mathbf{u}}_r$.

\begin{figure}[!t]
    \centering
    \includegraphics[width=.41\linewidth]{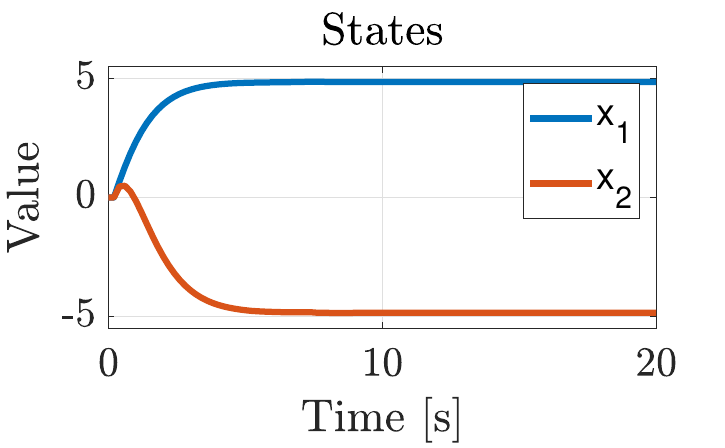}
    \includegraphics[width=.41\linewidth]{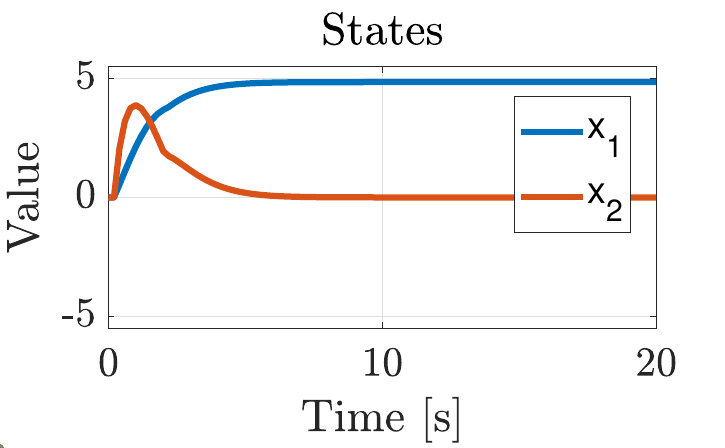}\\
    \includegraphics[width=.41\linewidth]{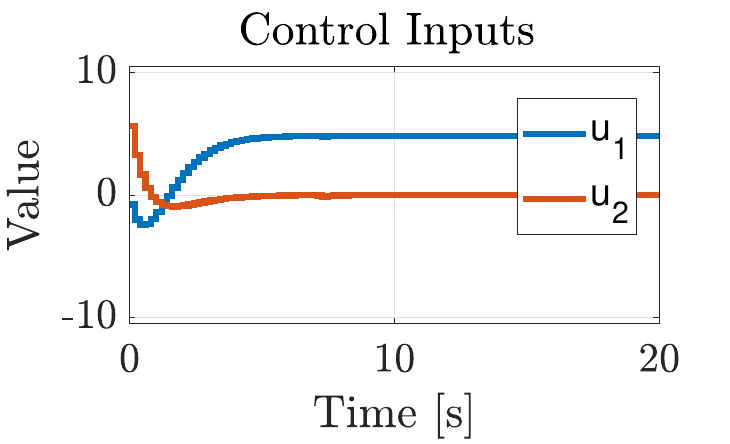}
    \includegraphics[width=.41\linewidth]{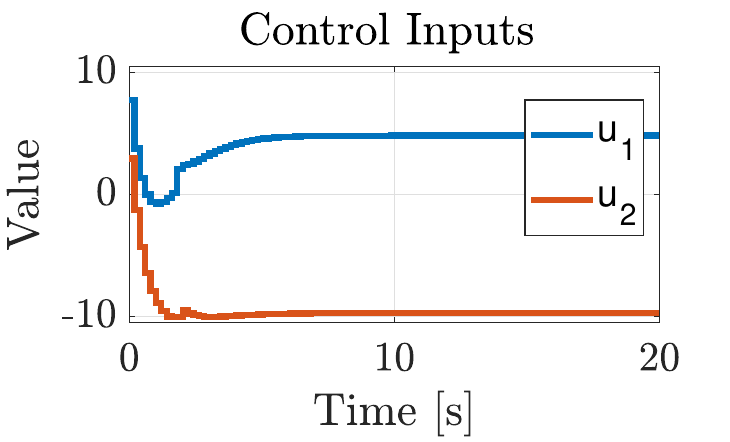}\\
    \includegraphics[width=.41\linewidth]{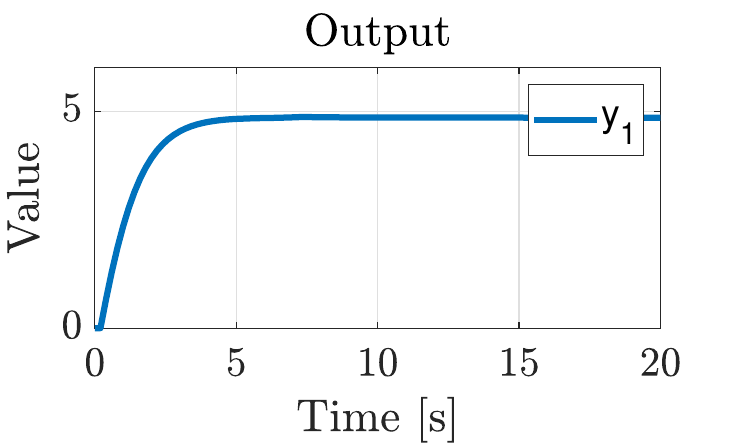}
    \includegraphics[width=.41\linewidth]{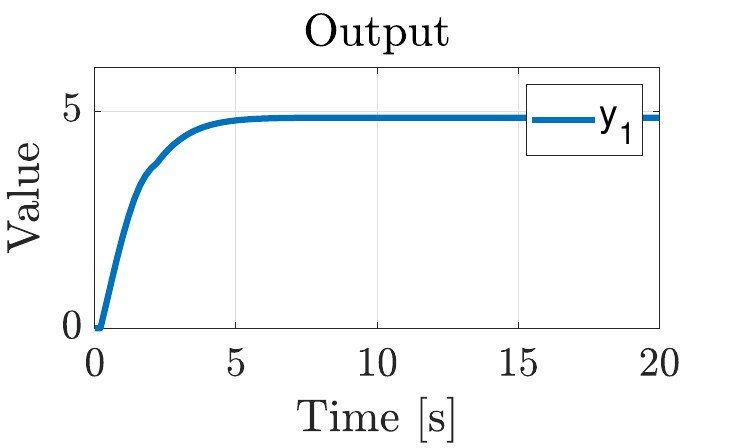}
    \caption{REAP-T Results for states, control inputs, and output of the system described in \citep{amiri2025practical} with $r=4.85$ (left), and desired steady state $\bar{\mathbf{x}}_r=[4.85~0]^\top$ (right).}
    \label{fig:TargetComparison}
\end{figure}

\begin{rem}
For a given desired reference $r$, there may be multiple associated steady-state configurations $\big(\bar{\mathbf{x}}_r, \bar{\mathbf{u}}_r\big)$. In such cases, REAP-T selects one steady-state configuration and guides the system toward it. To illustrate, consider the problem described in \citep{SSMPC} with the following scenarios: (1) the desired reference $r=4.85$ is specified; and (2) the desired steady-state configuration $\bar{\mathbf{x}}_r=[4.85~0]^\top$ is provided. As shown in Fig. \ref{fig:TargetComparison}, while the states and control inputs converge to different values in each case, the output consistently converges to the desired value of 4.85. 
\end{rem}

\begin{rem}
Since a core aspect of REAP theory involves tightening the constraints (see \citep{amiri2025practical} for more details), the desired steady state $\bar{\mathbf{x}}_r$ and steady input $\bar{\mathbf{u}}_r$ must lie strictly within the interior of $\mathcal{X}$ and $\mathcal{U}$, respectively. Therefore, users should avoid specifying steady-state configurations that are near the boundaries of the sets $\mathcal{X}$ and $\mathcal{U}$.
\end{rem}

    



\subsection{Terminal Constraint Set $\Omega$\textemdash Red Box in Fig. \ref{fig:GUI}}\label{sec:TerminalConstraintSet}
Let $\kappa(x(k), r) = \bar{\mathbf{u}}_r + K(x(k) - \bar{\mathbf{x}}_r)$ be a terminal control law, where $K = -(Q_u + B^\top Q_N B)^{-1}(B^\top Q_N A)$. The gain $K$ is designed to ensure that the closed-loop system matrix $A + BK$ is Schur. Consequently, the terminal constraint set $\Omega$, as defined in \eqref{eq:CostFunctionExtend}, is an invariant set under the terminal control law $\kappa(x(k),r)$ and is entirely contained within the system's constraints.

REAP-T employs two distinct approaches to implement the terminal constraint set $\Omega$, which will be discussed next, allowing users to select their preferred method via a dropdown menu.



\noindent\textbf{Prediction-Based:} If the pair ($C,A$) is observable, one approach to implement the terminal constraint set $\Omega$ is to use the maximal output admissible set, as described in \citep{Gilbert1991,SSMPC,amiri2025practical}. Specifically, the set $\Omega$ can be defined as:
\begin{align}\label{eq:MaxiamlOutputAdmissibleSet}
\Omega=\{(x,r):~\tilde{x}(\omega|x)\in\mathcal{X},\tilde{u}(\omega|x)\in\mathcal{U},~\omega=0,1,\cdots,\omega^\ast\},
\end{align}
where the sets $\mathcal{X}$ and $\mathcal{U}$ are defined as in \eqref{eq:allConstraints}, $\tilde{x}(0|x)=x$, $\tilde{x}(\omega|x)=(A + BK)^\omega x+\sum_{j=1}^\omega (A + BK)^{j-1} \left(B\bar{\mathbf{u}}_r - BK\bar{\mathbf{x}}_r\right)$ for $\omega\geq1$,  $\hat{u}(0|x)=K(x-\bar{\mathbf{x}}_r) + \bar{\mathbf{u}}_r$, and $\hat{u}(\omega|x)=K(A + BK)^\omega x+ K\sum_{j=1}^\omega (A + BK)^{j-1} \left(B\bar{\mathbf{u}}_r - BK\bar{\mathbf{x}}_r\right)+\left(B\bar{\mathbf{u}}_r - BK\bar{\mathbf{x}}_r\right)$ for $\omega\geq1$. In \eqref{eq:MaxiamlOutputAdmissibleSet}, $\omega^\ast$ is a finite index that can be computed by solving a series of offline mathematical programming problems as detailed in Algorithm \ref{al:Omegastar}. In REAP-T, the function \texttt{computeOmegastar} executes this algorithm to compute $\omega^\ast$.

\begin{algorithm}[!h]
\caption{Calculating the index $\omega^\ast$}\label{al:Omegastar}
\textbf{Input:} System matrices $A$, $B$, $C$, and $D$, terminal controller gain $K$, and steady-state configuration $(\bar{\mathbf{x}}_r,\bar{\mathbf{u}}_r)$.
\textbf{Notation Simplification:} For a given $\phi$, let constraints $\tilde{x}(\phi|x)\in\mathcal{X}$ and $\tilde{u}(\phi|x)\in\mathcal{U}$ be expressed as $2(n+p)$ inequalities of the form $\Phi_i(\phi|x)\leq0,~i\in\{1,\cdots,2(n+p)\}$.
\begin{algorithmic}[1]
\State Set $\phi=0$.
\State Solve the following optimization problems for $i=1,\cdots,2(n+p)$:
\begin{align*}
J_i^\ast=\left\{
\begin{array}{ll}
    & \max\Phi_i(\phi+1|x)\\
\text{s.t.} & \Phi_j(q|x)\leq0,~\forall j,~q=0,\cdots,\phi\\
\end{array}
\right..
\end{align*}

\State If $J_i^\ast\leq0,~\forall i$, set $\omega^\ast=\phi$. 

\State Otherwise, update $\phi\leftarrow\phi+1$ and return to Step 2.
\end{algorithmic}
\end{algorithm}

\begin{rem}
If the pair (C, A) is not observable, REAP-T displays an error message indicating that the user must change the method used to implement the terminal constraint, as follows:\vspace{0.1cm}
\begin{mdframed}
\texttt{The pair (C, A) is not observable. Please use the Lyapunov-based method to implement the\\ terminal constraint set.}
\end{mdframed}  
\end{rem}

\begin{rem}
REAP-T executes Algorithm \ref{al:Omegastar} until $\phi= 100$. Consequently, if REAP-T reports $\omega^\ast = 100$, it indicates that the prediction-based method may not be suitable for implementing the terminal constraint set for the specified system. In such cases, it is recommended to use the Lyapunov-based method, which will be discussed next.  
\end{rem}

\noindent\textbf{Lyapunov-Based:} An alternative approach for implementing the terminal constraint set $\Omega$ involves leveraging Lyapunov theory, as described in \citep{Hosseinzadeh2023RobustTermination,nicotra2018embedding}. Given that $A+BK$, with the specified $K$, is Schur, it follows that there exists a positive definite matrix $\Psi = \Psi^\top \succ 0$ ($\Psi \in \mathbb{R}^{n \times n}$) such that: $(A+B K)^\top \Psi(A+$ $B K)-\Psi \prec 0$. Using this property, the terminal constraint set $\Omega$ can be formulated as:
\begin{align}
\Omega=\big\{(x,r):~\left\Vert x-\bar{\mathbf{x}}_r\right\Vert^2_{\Psi}\leq\Gamma_i,i=1,\cdots,2(n+p)\big\},
\end{align}
where 
\begin{align}
\Gamma_i=\left\{
\begin{array}{rl}
   \frac{(\mathds{1}_i^\top\bar{\mathbf{x}}_r-\bar{x}_i)^2}{\mathds{1}_i^\top\Psi^{-1}\mathds{1}_i},  & i=1,\cdots,n \\
    \frac{(-\mathds{1}_i^\top\bar{\mathbf{x}}_r+\underline{x}_i)^2}{\mathds{1}_i^\top\Psi^{-1}\mathds{1}_i}, & i=n+1,\cdots,2n\\
    \frac{(\mathds{1}_i^\top K\bar{\mathbf{u}}_r-\bar{u}_i)^2}{\mathds{1}_i^\top K\Psi^{-1}K^\top\mathds{1}_i}, &  i=2n+1,\cdots,2n+p\\
    \frac{(-\mathds{1}_i^\top K\bar{\mathbf{u}}_r+\underline{u}_i)^2}{\mathds{1}_i^\top K\Psi^{-1}K^\top\mathds{1}_i}, &  i=2n+P+1,\cdots,2(n+p)
\end{array}
\right.,
\end{align}
with $\mathds{1}_i$ being a vector of appropriate dimensions, with its $i$-th element equal to 1 and all other elements equal to 0.

\textbf{Comparison:} The prediction-based method, which is limited to observable systems, implements the constraint \eqref{eq:ConstraintTerminalExtend} using $2(n+p)(\omega^\ast+1)$ linear constraints. In contrast, the Lyapunov-based method does not require observability and relies on only $2(n+p)$ quadratic constraints. This reduction in the number of constraints allows REAP-T to perform more iterations within a fixed time frame when the Lyapunov-based method is used. 

\begin{figure}[t]
    \centering
    \includegraphics[width=0.7\columnwidth]{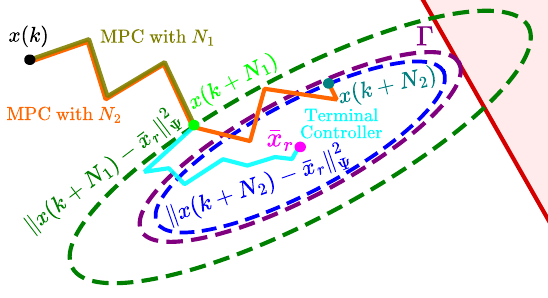}
    \caption{Illustration of prediction-based and Lyapunov-based methods for implementing the terminal constraint set, where $N_2>N_1$.}
    \label{fig:Terminals}
\end{figure}

However, the Lyapunov-based method tends to be more conservative as it relies on the Lyapunov level set at the end of the prediction horizon. Consequently, it often requires a longer prediction horizon to ensure its applicability, as illustrated in Fig. \ref{fig:Terminals}. 


\begin{rem}
If the specified prediction horizon length is insufficient for the Lyapunov-based method, REAP-T displays the following message and prompts users to increase the prediction horizon length:\vspace{0.1cm}
\begin{mdframed}
\texttt{The specified prediction horizon length is\\ insufficient for implementing the Lyapunov-\\
based method. Please increase the prediction\\
horizon length.}
\end{mdframed}
\end{rem}

\subsection{Reports and Plots\textemdash Violet Box in Fig. \ref{fig:GUI}}
\begin{figure}[t]
    \centering
    \includegraphics[width=0.4\linewidth]{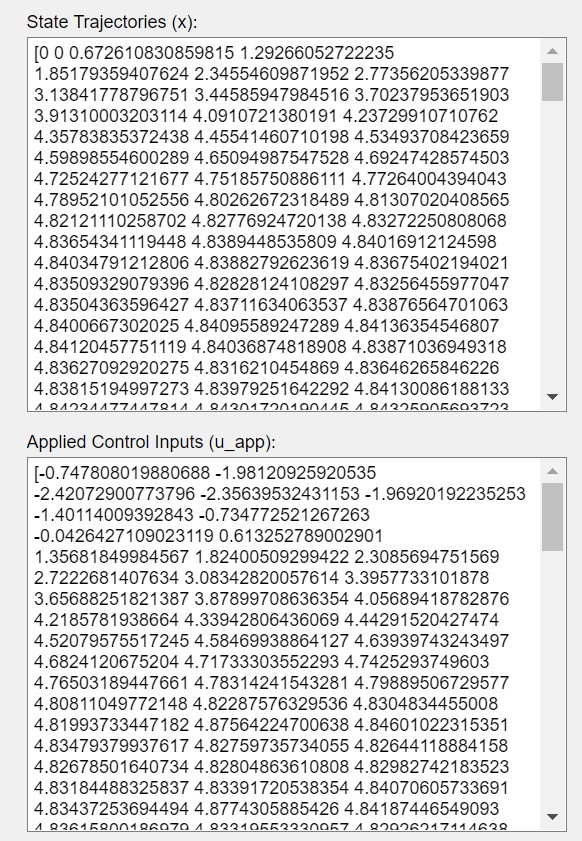}
    \caption{System states and inputs reported by REAP-T.}
    \label{fig:Reports}
\end{figure}

After the simulation, REAP-T provides a detailed report of system states and control inputs, as shown in Fig. \ref{fig:Reports}. Users can also generate plots of the results, including time profiles of the system states $x(t)$, control inputs $u$, system output $y(t)$, and the evolution of the KKT parameter $\sigma$ during REAP-T execution.



\subsection{Other Features and Functionalities}


\noindent\textbf{One-step Delay Implementation:}
REAP-T implements MPC under the Logical Execution Time (LET) paradigm \citep{Hosseinzadeh2022_ROTEC}, a widely adopted approach in cyber-physical systems. Specifically, REAP-T calculates the control signal using measurements taken at sampling instant $k$ and applies it to the plant at sampling instant $k+1$. This approach effectively implements MPC with a zero-order hold and a one-sample delay.










\noindent\textbf{Warm Starting:} At each time step, REAP-T utilizes the warm-starting strategy outlined in \citep{Hosseinzadeh2023RobustTermination}, implemented via the \texttt{Warmstarting} function, to initialize both the primal and dual variables.

For the primal variable, warm starting involves a one-step backward shift of the previous control input, combined with the terminal control law, to determine values for the new prediction instant. For the dual variable, REAP-T similarly uses the dual variables from the previous time instant as initial values for the new prediction.

\noindent\textbf{Acceptance/Rejection Mechanism:} As discussed in \citep{Hosseinzadeh2023RobustTermination} and \citep{amiri2025practical}, the evolution of $\big(\hat{\mathbf{u}}(\tau|k) - \mathbf{u}^{\ast}(k)\big)$ may not be monotonic. Consequently, early termination of the virtual dynamical system \eqref{eq:DisREAP} could compromise the closed-loop stability of the system. To address this issue, REAP-T employs a logic-based acceptance/rejection mechanism implemented through the function \texttt{ARMechanism}.

At each time instant $k$, REAP-T begins with the initial condition $\big(\hat{\mathbf{u}}(0|k),\hat{\lambda}(0|k)\big)$ and executes $\tau_k$ computation iterations. The control input $\hat{\mathbf{u}}(\tau_k|k)$ is accepted if the value of the cost function defined in \eqref{eq:CostFunctionExtend} is lower than that obtained with $\hat{\mathbf{u}}(0|k)$.

\section{Tutorial}\label{sec:Tutorial}
This section provides a brief overview of two examples included in REAP-T.

\subsection{Position Control of the Parrot Bebop 2 Drone}
\begin{figure}[!t]
    \centering
    \includegraphics[width=0.5\linewidth]{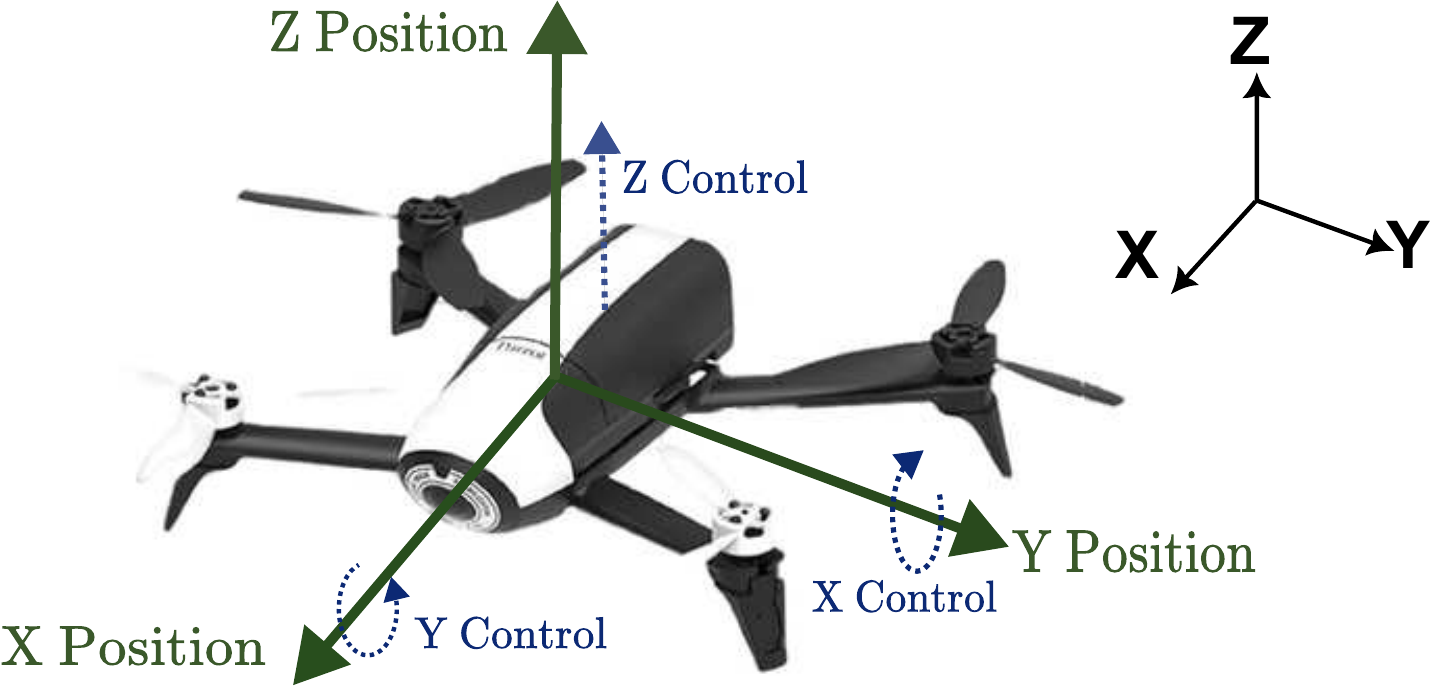}\vspace{.5em}\\\includegraphics[width=0.41\linewidth]{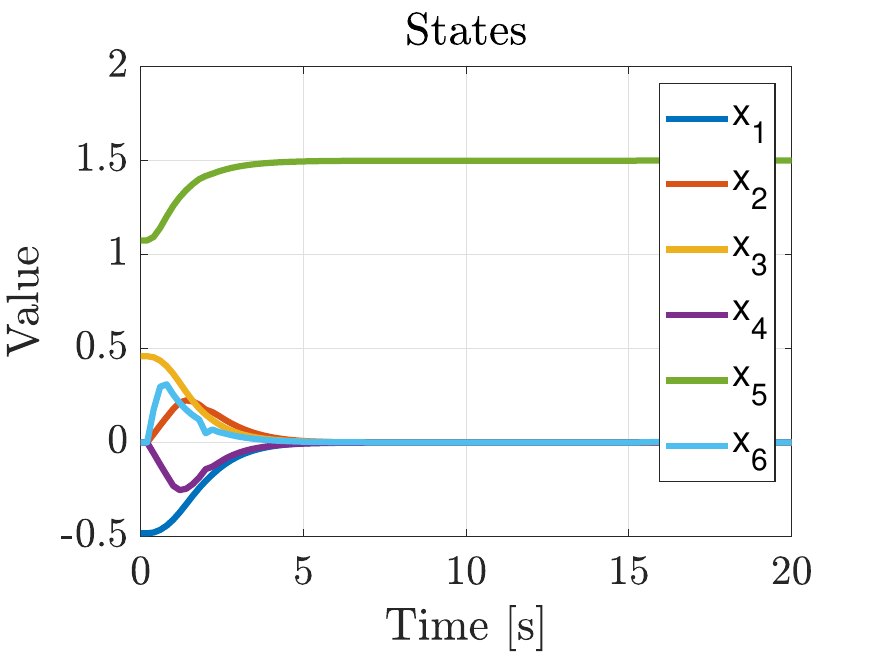}\includegraphics[width=0.37\linewidth]{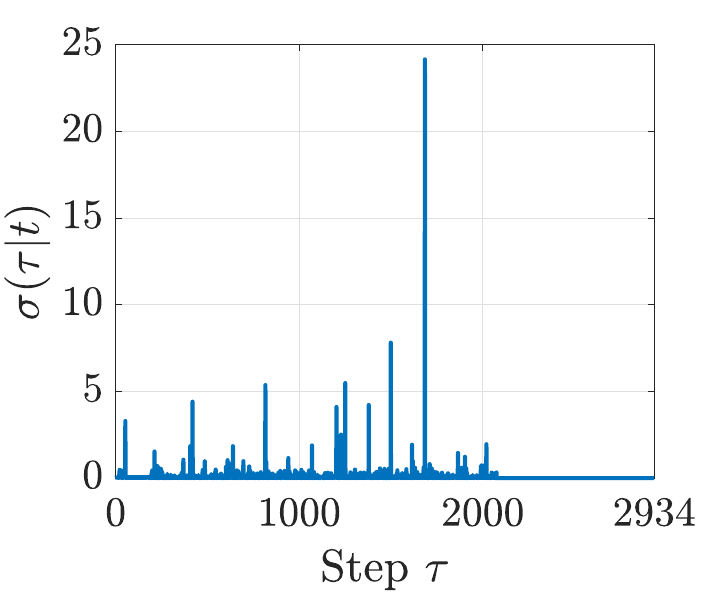}
    \caption{Parrot Bebop 2 with considered reference and body frames, and REAP-T results for the position control.}
    \label{fig:Output}
\end{figure}

The first example involves a Parrot Bebop 2 drone, whose continuous-time dynamical model is described in \citep{AMIRI2024330}. The system includes six states and three control inputs, with the following constraints: X constraint U.B.$=[10~10~2.57~10~10~10]^\top$, X constraint L.B.$=[-10~-10~-10~-10~0~-10]^\top$, U constraint U.B.$=[0.05~0.05~0.6]^\top$, and U constraint L.B.$=[-0.05~-0.05~-0.6]^\top$. We discretize the system with a sampling period of $\Delta T=0.2$, and set the weighting matrices to $Q_x = \text{diag}\{5\times I_4, 1000\times I_2\}$,$Q_u = \text{diag}\{35, 20, 1\}$. Fig. \ref{fig:Output} presents the results with initial condition $x(0)=[-0.48~0~0.46~0~1.08~0]^\top$
and the reference equilibrium point $\bar{x}_r = [0~0~0~0~1.5~0]^\top$.

\subsection{Roll and Side-slip Angles Control of F-16 Aircraft}
\begin{figure}[!t]
    \centering
    \includegraphics[width=0.5\linewidth]{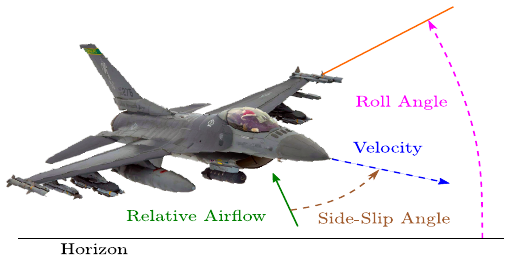}\\\includegraphics[width=.41\linewidth]{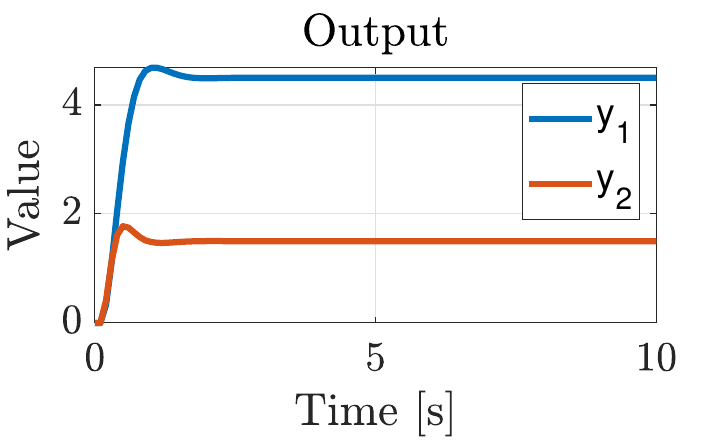}
    \includegraphics[width=.41\linewidth]{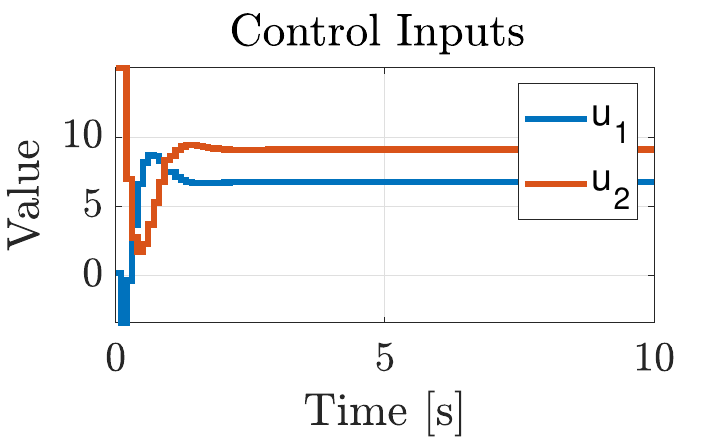}
    \caption{F-16 aircraft, and REAP-T results for controlling the roll and slide-slip angles.}
    \label{fig:F16}
\end{figure}

The second example features a simplified continuous-time model representing the lateral dynamics of an F-16 aircraft, as described in \citep{suresh2005nonlinear}. This system consists of four states and two control inputs, which are subject to the following constraints: X constraint U.B.$=[\text{Inf}~\text{Inf}~5~2]^\top$, X constraint L.B.$=[-\text{Inf}~-\text{Inf}~-5~-2]^\top$, U constraint U.B.$=[10~15]^\top$, and U constraint L.B.$=[-10~-15]^\top$. We discretize the system with a sampling period of $\Delta T=0.1$, and set the weighting matrices to $Q_x = \text{diag}\{0.1,0.1,10,10\}$,$Q_u = \text{diag}\{0.1,0.1\}$. Figure \ref{fig:F16} presents the output and control inputs generated by REAP-T when the desired reference is set to $r=[4.5~1.5]^\top$.





\section{Conclusion}\label{sec:Conclusion}
This paper provided a brief overview of the theoretical framework for robust-to-early termination MPC and introduced a MATLAB toolbox, which is named REAP-T, for its implementation. REAP-T is modular, user-friendly, and designed to allow users to explore, analyze, and customize its features and functionalities. Additionally, two real-world examples are included to guide users on effectively utilizing the toolbox. Future versions of the toolbox will incorporate additional examples, algorithms, and features to expand REAP's applicability and versatility.

\section*{Acknowledgment}\vspace{-0.3cm}
This research has been supported by National Science Foundation under award number DGE-2244082.


\bibliography{ref}        

\end{document}